\documentclass{article}
\usepackage{ amssymb }
\usepackage[utf8]{inputenc}
\usepackage{lipsum}
\usepackage{balance}
\usepackage{graphicx}
\usepackage{nccmath}
\usepackage{subfig}
\usepackage{mathtools}
\usepackage{times}
\usepackage{url} 
\usepackage{appendix}
\usepackage{caption}
\usepackage{amsmath}
\usepackage{amsthm}
\newtheorem{theorem}{Theorem}[section]

\newtheorem{lemma}[theorem]{Lemma}
\usepackage{algorithm}
\usepackage{algpseudocode}
\usepackage{float}
\floatstyle{plaintop}
\usepackage[margin=1in]{geometry}

\title{Preserve Prime}

\begin{document}
\title{Prime Holdout Problems}
\author{Max Milkert, Alex Ruchti}
\maketitle
\begin{abstract}
The Collatz conjecture, otherwise known as the $3n +1$ problem, is simple to state, yet has gone unsolved for over 70 years. Many researchers have attempted to reason about the conjecture by looking at generalizations, hoping to have results that offer a solution to the specific problem. This paper takes a different approach, introducing a class of problems opposite the Collatz conjecture that we refer to as holdout problems. The difference is that, after applying a linear function, instead of dividing by a finite set of prime factors, a holdout problem specifies a set of primes to be retained. A proof that all positive integer starting values converge to 1 is given for an example of a finite and an infinite holdout problem. It is conjectured that finite holdout problems cannot diverge on any input, which has implications for divergent sequences in the Collatz conjecture.
\end{abstract}

\section{Introduction}
The Collatz Conjecture, also known as the 3n + 1 problem, is a simply-formulated problem that is notoriously unproven. It concerns the iteration of a piecewise equation defined on the set of positive integers \cite{lagarias2010ultimate}

\begin{ceqn}
\begin{align}
    f(n) = \begin{cases} 
      3n + 1 & n \equiv 1 \pmod{2} \\
      \frac{n}{2^{v_2(n)}} & n \equiv 0 \pmod{2} \\
   \end{cases}
\end{align}
\end{ceqn}
where $v_2(n)$ is the maximum power of 2 dividing $n$ (also called the 2-adic order of $n$).\\*

To better illustrate the behavior of this function, consider the sequence formed by applying $f(n)$ beginning at 7: $7 \rightarrow 22 \rightarrow 11 \rightarrow 34 \rightarrow 17 \rightarrow 52 \rightarrow 13 \rightarrow 40 \rightarrow 5 \rightarrow 16 \rightarrow 1 \rightarrow 4 \rightarrow 1... $ We see that repeated applications of $f(n)$ alternate between removing all factors of 2 from $n$ and applying the 3$n$ + 1 linear map, eventually entering into a cycle between 1 and 4, often called the trivial cycle.\\*

Lothar Collatz famously postulated that for all starting values $n$, iteration of $f(n)$ will always reach 1 (or equivalently descend into the trivial cycle). And while this conjecture has been circulating since the 1950s \cite{lagarias2010ultimate}, all attempts at proving it have been unsuccessful. A successful proof would have to do two things: rule out the possibility that there could be non-trivial cycles, and rule out the possibility that there is a starting value for which the iteration would diverge. Either of these could also be an avenue for disproof, in the form of a cycle this could be fairly straightforward, one would only need to find the cycle to have a valid counterexample. Showing that a sequence is divergent, on the other hand, could be a more complicated avenue of disproof. A counterexample of this nature would be infinitely long, thereby requiring a full proof to reason about 1's absence. In light of current open nature of the problem, despite many attempts at a solution, there's also a reasonable suspicion that the problem might not be decidable at all\\*

\subsection{Periodic Linear Functions}
Some authors have chosen to focus on generalizations of the problem, in the hopes that a finding about a larger problem class can be applied back to the Collatz conjecture itself. John H. Conway (known for Conway's game of life) has examined the problem as an instance of a periodic linear function \cite{Conway} to argue for his suspicion that there might not be a solution. A periodic linear function is a function of the form:

\begin{ceqn}
\begin{align}
    P(n) = \begin{cases} 
      a_0n + b_0 & n \equiv 0 \pmod{N} \\
      a_1n + b_1 & n \equiv 1 \pmod{N} \\
      \vdots\\
      a_{N-1}n + b_{N-1} & n \equiv N-1 \pmod{N}
   \end{cases}
\end{align}
\end{ceqn}
where $N$ is some arbitrary modulus, and the $a$ and $b$ coefficients are rational numbers chosen so that $P(n)$ will always be an integer. The term "periodically linear" arises from the fact that each of these component functions $a_in + b_i$ is linear, and the decision on which to apply repeats periodically over the number line.\\

The Collatz function may be encoded in this form as
$N = 2$, $a_0 = 1/2$, $b_0 = 0$, $a_1 = 3$, and $b_1 = 1$. Conway was able to show that each member of this class with all $b_i =0$ can be converted into a Minsky machine (and vice versa). A Minsky machine can be thought of as a computer program, it comes equipped with memory variables capable of holding any positive integer, a set of execution states, and rules that govern how it transitions between those states to modify the memory with the goal of computing the result of some problem. Minksy machines bear semblance to the more well known Turing machines, and are proven to have equal computational power. This means that many of the major results from computability theory can carry over to periodic linear functions. Of chief concern in this paper is the general undecidability of the halting problem. It is not possible to produce a program that can correctly state for all programs whether or not they terminate on a given input. The halting problem is strongly related to the Collatz conjecture since we could consider each sequence to simply terminate upon reaching 1, in which case the questions of whether each starting value returns to 1, and whether the Collatz iteration halts for every input are identical. Since the set of all periodic linear functions contains the set of all programs, then there have to be examples of problems in this class for which the halting problem (and thus whether or not 1 is always reached) is undecidable.\\*

\subsection{Divisor Problems}

While this may seem like bad news for the Collatz conjecture, recall that Conway's argument shows a bijection with the Minsky machines only between those periodic linear functions with all $b_i = 0$. While the collatz function is periodically linear, is not a member of this narrower class, since it has $b_1 = 1$. Therefore it might make sense to consider an alternate generalization of the conjecture, ideally this generalization would be disjoint from the $b=0$ set. One possible way of doing this would be to give the Collatz function the freedom to use a number other than 3, and to allow a set of prime divisors $\mathbb{D}$ that is potentially larger than just $\{2\}$. Throughout the paper we will refer to these as divisor problems, let 

\begin{ceqn}
\begin{align}
    C_{0}(n) = \frac{n}{d(n)}
\end{align}
\begin{align}
    C(n) = C_{1}(n) = \frac{(cn+1)}{d(cn+1)}
\end{align}
\begin{align}
    C_{k+1}(n) = C(C_{k}(n))
\end{align}
\end{ceqn}
where $c \in \mathbb{N}$ is a constant and
\begin{ceqn}
\begin{align*}
d(n) = \prod_{d_i \in \mathbb{D}} {d_i}^{v_{d_i}(n)}
\end{align*}
\end{ceqn}
$v_{d_i}$ represents the $d_i$-adic order of $n$ for divisor $d_i$.\\*

A few notes on this definition: the result of $d$ captures all needed division in one step, this definition of $C$ enables us to avoid using a piecewise definition as in (1). $C_0$ will provide an initial removal of all prime factors in $\mathbb{D}$, and all further applications of $C$ will keep prime factors in $\mathbb{D}$ out of the result. Like the Collatz conjecture, a divisor problem asks whether or not for all integer starting values $n$, there exists some $k$ such that $C_{k}(n)=1$.\\*

There are examples of divisor problems being studied in the Collatz literature. Recalling from earlier, one of the 2 ways to show the falsehood of the Collatz conjecture would be to produce a non-trivial cycle. In his search for such a cycle, Tomas Oliveira e Silva developed an efficient computer assisted approach to determining their existence \cite{eSilva}. In that paper he also extended his approach to also detect if there are cycles in other divisor problems. One such problem is $C_{5,\{2,3\}}$, referring to the $5n + 1$ divisor problem in which the divisor set is $\{2,3\}$, a notation we will adopt for use in this paper. Computational approaches to the Collatz conjecture continue to this day \cite{roosendaal}, the reader may even lend some of their CPU (or GPU) cycles to the cause should they be so inclined. At the time of the writing all starting values below $2^{69}$ have been determined to converge.\\*

\subsection{Holdout Problems}

There are many other possible expansions of the Collatz conjecture than just those discussed here, each aiming to better characterize certain features of the problem. This paper will take a slightly different approach, to attack the Collatz problem not through a generalization, but rather by examining its opposite. Here we introduce holdout problems, so called because in each iteration, a finite set $\mathbb{H}$ of prime factors are held aside while the rest are discarded, these numbers we will refer to as a holdout set. For a given set of parameters, a holdout and a divisor problem will retain exactly the opposite set of factors. This connection between the two classes is evident from their definitions. Let
\begin{ceqn}
\begin{align}
G_0(n)= h(n)
\end{align}
\begin{align}
    G(n) = G_{1}(n) = h(cn+1) 
\end{align}
\begin{align}
    G_{k+1}(n) = G(G_{k}(n))
\end{align}
\end{ceqn}
where $c\in \mathbb{N}$ is again constant and
\begin{ceqn}
\begin{align*}
\text{ }h(n) = \prod_{h_i \in \mathbb{H}} {h_i}^{v_{h_i}(n)}
\end{align*}
\end{ceqn}

You will notice the function $h$ is identical to $d$ of the divisor problems, with the difference between the two problem classes being in how the result is used. Similarly, $G$ both performs a linear mapping, and incorporates an application of $h$ into each iteration. Like a divisor problem, a holdout problem also asks whether or not for all integer starting values $n$, there exists some $k$ such that $G_{k}(n)=1$\\*

To see a holdout problem in action, we can consider the sequence of $H_{5,\{2,3\}}$ beginning at $n=28$. $28 = 4\times7$, 4 is retained since it is a power of 2, 7 is discarded since it is not a member of $\{2,3\}$. 4 will undergo the linear map $5\times4 + 1 = 21$. $21 = 3\times7$, as a power of 3, 3 is retained and 7 is again discarded. This process continues alternating until 1 is reached. The sequence in its entirety is as follows: $28 \rightarrow 4\rightarrow 21 \rightarrow 3 \rightarrow 16 \rightarrow 81 \rightarrow 406 \rightarrow 2 \rightarrow 11 \rightarrow 1$. From the definition, this sequence should only contain terms divisible by 2 or 3, but we are including intermediate terms before applying $h()$ to better illustrate the workings of this process. It can also be noted that sometimes $h(n) = 5n +1$ and thus the intermediate terms are absent, like the case for $3 \rightarrow 16$.\\*

$H_{5,\{2,3\}}$ is the primary problem that we will be studying in this paper, and a proof of its convergence will be presented in section 2. This problem was selected instead of the "opposite" of the Collatz conjecture. This is because it can be readily observed that $H_{3,\{2\}}$ converges: odd numbers, lacking any factor of 2, will immediately be sent to 1. Even numbers will become powers of 2, which are mapped to odd numbers by $3n + 1$. To try to get more complex behavior, we could try expanding the holdout set to include 3 ($H_{3,\{2\}}$) but $3n+1$ is never divisible by 3, so after the first iteration, it's back to following $H_{3,\{2\}}$. It is only once we move to $H_{5,\{2,3\}}$ that we begin to see interesting behavior like the sequence above.\\*

\subsection{Infinite Holdout Problems}

In between the holdout and divisor problems lies a third class of problems, the infinite holdout problems. This is a class distinct from divisor problems, even though a divisor problem can be translated into a holdout problem with an infinite holdout set. The distinguishing factor is that in an infinite holdout problem, both the holdout and divisor sets are infinite. An example of an infinite holdout problem is the holdout problem where
\begin{ceqn}
\begin{align}
    \text{ } c = 2 \text{ and }\mathbb{H} = \{p \in \mathbb{P} \mid p \equiv 1\pmod{4}\}
\end{align}
\end{ceqn}
where $\mathbb{P}$ denotes the set of all primes\\*

It is famous historical result in number theory that both the set of primes congruent to 1 and congruent to 3 mod 4 are infinite. Although this may appear more complicated than a finite holdout problem, section 3 presents a surprisingly short proof of its convergence. This paper concludes with section 4, where we further discuss the relation of holdout and divisor problems to Conway's problem class and the theory of computation, and explore how holdout problems might be employed to produce results about the Collatz conjecture.\\*

\subsection{Terminology}
Throughout this paper we make use of several terms and symbols. The terms "orbit" and "trajectory" can sometimes be used to describe sequences in the Collatz conjecture, and the individual terms can be referred to as iterates. Cycling will by default describe non-trivial cycles, and convergence will mean convergence to 1 unless otherwise stated.  $\mathbb{H}$ and $\mathbb{D}$ will be used to refer to the holdout and divisor sets of a problem. The notation $H_{c,\mathbb{H}}$ will refer to the holdout problem that applies $cn+1$ and preserves the set $\mathbb{H}$ between iterations.

\section{Proof of convergence for $H_{5,\{2,3\}}$}

\begin{lemma}{if $a \mid n$, $a \nmid G(n)$}
\label{add_one}
\begin{proof}
The result is immediate, since $a \mid n$ we have
\begin{gather*}
    5 \times n + 1 \equiv 1 \pmod{a}
\end{gather*}
As a direct consequence of this lemma, $G(6n) = 1$, as it cannot be divisible by $2$ or $3$.\\*
\end{proof}
\end{lemma}

It can be observed that all sequences of substantial length consist of alternating powers of either prime number in the holdout set. Stated precisely it gives the following lemma:

\begin{lemma}{if $G_k(n) \neq 1$ for all $0 \leq k \leq K$ and $K > 1$, then the iterates $G_0(n)...G_{K}(n)$ are composed of alternating powers of $2$ and $3$}
\label{alternating}
\begin{proof}
The restriction that $1$ does not appear in the sequence (until possibly after the $K$th term) allows us to exclude several cases from consideration. Sequences with starting values not divisible by $2$ or $3$ are excluded since $G_0(n) = h(n) = 1$. Additionally, by Lemma \eqref{add_one}  $G(6n) = 1$ so none of the iterates $G_0(n)...G_{K-1}(n)$ are divisible by $6$\\

The only remaining possible factorization for the $G_k(n)$ are $2^a$ or $3^b$. The alternation between the powers follows from Lemma \eqref{add_one} since $3 \nmid G(3^a)$ and $2 \nmid G(2^b)$. By the same reasoning the final term $G_K(n)$ cannot be divisible by $6$
\end{proof}
\end{lemma}

\begin{lemma}{ For all $a \ge 0$, $G(2^{2a+1}) = 1$ }
\label{Odd2}
\begin{proof}
\begin{gather*}
    (5 \times 2^{2a + 1}) + 1 \equiv (2\times 1 \times 2) + 1 \not\equiv 0 \pmod{3}.
\end{gather*}
Since $2 \nmid (5 \times 2^{2a + 1} + 1)$ we must have $G(2^{2a + 1}) = 1$. Therefore all trajectories reaching an odd power of 2 are convergent.\\* 
\end{proof}
\end{lemma}

A rather unique property of $\mathbb{H}_{5,{2,3}}$ is that when $G(n)$ produces a larger number than $n$, then no division may have occurred in that step, which is expressed in the following lemma: 
\begin{lemma}{ for all $n\ge 1$ we have either $G(n)<n$ or $G(n) = 5n+1$}
\label{DivisionLemma}
\begin{proof}
Assume $G(n)\neq 5n+1$. Let $5n+1 = 2^a 3^b c$ where $a = v_2(5n + 1)$ and $b = v_3(5n + 1)$. This ensures that $2,3 \nmid c$. Applying the function $h$ will remove $c$ from this factorization.\\*

Using the above we can rule out all values of $c$ that are less than $7$. $G(n)\neq 5n+1$, so $c \neq 1$. We have $2 \nmid c$ and $3 \nmid c$, eliminating $2,3,4,$ and $6$ as possibilities. $5 \nmid c$ otherwise $5n +1 = 2^a 3^b c \not\equiv 1 \pmod{5}$. Therefore we have $c \ge 7$.
\begin{gather*}
     G(n) = \frac{5n + 1}{c} < n \text{ } for \text{ } all \text{ } c \ge 7 \text{ } and \text{ } n \ge 1\\
\end{gather*}
Therefore any time $n$ increases, it must be to $5n +1$ and no other value.
\end{proof}
\end{lemma}

The following 2 Lemmas will make use of this to show that there are only a few unique instances of $n$ increasing.\\*

\begin{lemma}{$G(3^a) < 3^a$ for all $a > 1$}
\label{Power3}
\begin{proof}
Assume $G(3^a) > 3^a$. By Lemma \eqref{DivisionLemma} we have $G(3^a) = 5 \times 3^a + 1$. By Lemma \eqref{alternating} we have that $G(3^a) = 2^B$.
Examining the equivalence $5 \times 3^a + 1 = 2^B$ modulo $3$ reveals that $2^B \equiv 1 \pmod{3}$ therefore $B$ is even. Let $B = 2b$, rearranging we have 
\begin{align*}
    5 \times 3^a = 2^{2b} - 1 
\end{align*}
factorizing this difference of squares gives us 
\begin{align*}
    5 \times 3^a = (2^b + 1) \times (2^b - 1) 
\end{align*}
Notice that only one of these terms may be divisible by $3$ because their difference is $2$, and that $5$ only appears once on the left side. By inspection we can see that the only solution is $a = 1$ and $b = 2$ when we assign
\begin{align*}
    5 = (2^b + 1) \text{ and } 3 = (2^b - 1) 
\end{align*}
Using any larger power of $3$ will make the difference larger than $2$, while having both $3$ and $5$ divide the same term will not produce any solutions.\\*
\end{proof}
\end{lemma}

\begin{lemma}{$G(4^a) < 4^a$ for all $a \neq 2$}
\label{Power4}
\begin{proof}
Assume $G(4^a) > 4^a$. By Lemma \eqref{DivisionLemma} we have $G(4^a) = 5 \times 4^a + 1$. By Lemma \eqref{alternating} we have that $G(4^a) = 3^B$.
Examining the equivalence $5 \times 4^a + 1 = 3^B$ modulo $4$ reveals that $3^B \equiv 1 \pmod{4}$ therefore $B$ is even. Let $B = 2b$ rearranging we have 
\begin{align*}
    5 \times 4^a = 3^{2b} - 1 
\end{align*}
factorizing this difference of squares gives us 
\begin{align*}
    5 \times 4^a = (3^b + 1) \times (3^b - 1) 
\end{align*}
Notice that the difference between the terms is again $2$, since the product is even, we can conclude that one will have remainder $2$ mod $4$ and the other will be divisible by $4$ (remainder $0$). Since their product is nonzero we must have $a \ge 2$ to produce all the required copies of $2$. Upon inspection we can see the solution is $a = 2$ and $b = 2$ assigning
\begin{align*}
    5 \times 2 = (3^b + 1) \text{ and } 2^3 = (3^b - 1) 
\end{align*}
again a larger power of $2$ will make the difference too large, as will assigning $5$ to the other term, this can be the only solution.\\*
\end{proof}
\end{lemma}

\begin{theorem}{$H_{5,\{2,3\}}$ converges for every integer starting value}
\begin{proof}
We begin by observing that Lemmas \eqref{add_one} and \eqref{alternating} allow us to reduce the scope of our proof to only powers of $2$ and $3$ since other starting points converge. Lemma \eqref{Odd2} allows us to further refine our proof to powers of $3$ and $4$.\\

We will proceed by induction using the sequence $4 \rightarrow 3 \rightarrow 16 \rightarrow 81 \rightarrow 2 \rightarrow 1$ as our base case. Fortunately this sequence demonstrates the convergence of all starting points less than or equal to the solutions in Lemmas \eqref{Power3} and \eqref{Power4}, strengthening them to ensure $G(n) < n$ outside that base sequence. We can now proceed using the induction hypothesis that $H_{5,\{2,3\}}$ converges for all numbers less than $n$ and our result follows.
\end{proof}
\end{theorem}

\section{Proof of convergence for $H_{2,\{p|p \equiv 1 \pmod{4}\}}$}

Here we show that the infinite holdout problem which retains all prime factors congruent to 1 modulo 4 converges to 1 for all starting values\\*

\begin{theorem}{The 2n+1 infinite holdout problem converges}
\begin{proof}
We proceed by induction, our base case being that 1 = 1. Let $n$ be a positive integer greater than 1. Assuming that all previous numbers converge to 1, there are 2 cases to consider:
\begin{enumerate}
    \item{$n$ is divisible by a prime $p$ congruent to 2 or 3 mod 4.}
    Thus $n$ can be divided, and converges because $n/p$ converges.
    \item{$n$ is composed only of prime factors congruent to 1 mod 4.}
    Since $n$ is odd, $2n+1$ is congruent to 3 mod 4. It must have at least one prime factor that is congruent to 3, the smallest of which is 3. For all $n$ greater than 1, $(2n+1)/3$ is less than $n$ therefore $n$ converges.
\end{enumerate}
\end{proof}
\end{theorem}

\section{Discussion}

In the previous sections we have demonstrated the truth of both a finite and an infinite holdout problem. In this section we examine holdout problems in light of Conway's results, and argue for why divisor problems might be a more promising way to view the Collatz conjecture. Lastly, but perhaps most significantly, we posit the finite holdout conjecture- that no finite holdout problem may support divergent trajectories- and discuss its implications for the Collatz conjecture.

\subsection{Conway's Theorem and Divisor Problems}

Recall from section 1.1 that Conway has shown the equivalence between Minsky machines and addition-free periodic linear functions (all $b_i = 0$). A key idea behind Conway's proof is that a Minsky machine stores whole numbers, rather than binary variables, in its memory. This property will allow for the execution state of a Minsky machine to be embedded as a single positive integer on which a periodic linear function may be iterated. Representing this number as the powers present in its prime factorization (i.e. a,b,and c... in $2^a,3^b,5^c...$) we may assign each prime power to hold the contents of one of the Minsky machine's memory registers. We can then set aside one additional power to track which state the machine is in. Once we have done this, we can see that the rational coefficients $a_{i}$ from the piecewise definition of $P(n)$ are acting as program instructions each time the function is applied, incrementing and decrementing the registers, and transitioning between states. A modulus $N$ can then be chosen such the remainders modulo $N$ provide information about the contents of the registers (perhaps by setting $N$ to be the product of the primes) so that valid transition rules can be defined.\\*

With this framework in mind, Conway was able to produce a periodic linear function using 24 different fractions that is universal \cite{fractran}. A universal Minsky machine, can be thought of as a full fledged computer, rather than just a single program. It is able to be given another Minsky machine and perfectly simulate its behavior by systematically stepping through its internal state machine and executing its logic. The notion of encoding the logic of a Minsky machine in an integer is rather strange, since we've only seen how to encode the memory thus far, but in \cite{fractran} Conway shows how it may be done, although it is rather complicated and the numbers involved typically end up being very large.\\*

Universal programs are of great theoretical importance, and play an important role in reasoning about the halting problem. This is because they can be used to craft an adversarial program against any possible solution attempt, rather than just one. Note that a solution to the halting problem need not be universal itself, as fully simulating a non-halting program would be futile. The solution does, however, have to take a program as a parameter and reason about its logic to produce a decision. Now to craft a specific adversarial program, we could take this solution and modify it as such: after executing the halting function on itself, the adversary could then choose to enter an infinite loop, or terminate such that the original halting program always reasons incorrectly. This would mean the halting problem is undecidable in general, since no program can provide the correct answer for whether its adversary terminates. A universal adversary would be even more sophisticated, able to simulate any halting function given to it as a parameter, it could defy that program's predictions about it so that even for this one program alone, determining if it halts on every input is impossible. This is all the more remarkable bearing in mind that this complicated chain of passing programs as parameters collapses into a single innocent seeming periodic linear function, started on some large number.\\*

The Collatz function itself, fortunately, is thought not to be universal \cite{lagarias2010ultimate}. Using Conway's construction, we will argue that it is unlikely that any other divisor problem is either. First we must convert a divisor problem into the form of a periodic linear function. Let

\begin{align*}
     N = \prod_{d \in \mathbb{D}}
\end{align*}
and let
\begin{align*}
     P(n) = \begin{cases} 
      cn + 1 & gcd(i,N) = 1 \\
      \frac{1}{gcd(i,N)}n + 0 & gcd(i,N) > 1 
   \end{cases}
\end{align*}
where $n \equiv i \pmod{N}$.\\ 
This periodic linear function will take the $cn + 1$ transition whenever the original divisor problem would have, but may have to take multiple steps to perform all the necessary division in between.\\*

Now suppose we were to encode a Minsky machine in a similar fashion to Conway, our program registers would then be tracking the powers of primes in $\mathbb{D}$. This machine has two notable features that would seem to make computation neigh impossible. The first is that there is a transition that adds 1. Were this to occur in a Minsky machine, it would be problematic. Any register (including the execution state) with a nonzero value would be wiped to 0. Additionally, any registers that were empty have the potential to be flooded with unknown contents. These contents are unknown since the program logic is dictated by a congruence modulo $N$, causing the program to behave identically regardless of the powers on the primes it is not observing. The second problem with divisor computation lies in all the other transitions, which dictate that the contents of the registers which the program tracks should all be zeroed out, which would also make computing anything meaningful extraordinarily difficult. Putting these two defects together, a divisor problem will repeatedly zero its memory, and then bring in random new memory values, before zeroing them again. Some divisor problems (the ones that converge) are extremely potent memory erasers, zeroing all prime powers, even the ones they were not observing. In light of this, it's hard to see how a universal model of computation could be encoded as a divisor problem, and if the intuition here is correct, then it closes one avenue by which the halting problem may be shown to be undecidable for this class. It could always be the case that computation is somehow embedded in the structure of the number line, that the effect of adding 1 on the divisibility by each prime has some underlying computational nature that could be exploited by carefully choosing a $c$ and a divisor set, but this would be a rather remarkable finding. It is also possible that there is some other reason that the convergence of divisor problems in general is undecidable.  

\subsection{Holdout Problems and Computability}

We have now seen why it is unlikely that Conway's results can carry over to divisor problems. Here we will see that the outlook is even brighter for holdout problems, because unlike divisor problems, neither finite nor infinite holdout problems are periodically linear. Both of these classes have infinite sets of divisors, trying to put them in a periodically linear form would require finding some finite modulus $N$ under which all their operations can be expressed. But this is not possible with infinite divisor sets, for any $N$ chosen, there will always be a prime $p$ in $\mathbb{D}$ that is larger than $N$. This periodic linear function will either fail to divide by $p$ when it is possible, or it will divide by $p$ when the result would not be an integer, since both numbers divisible by $p$ and not can end up with the same remainder modulo $N$.\\*

The sets of all finite holdout problems and all divisor problems are both countably infinite. For every prime in the holdout set (or the divisor set respectively), a binary variable could be assigned, each value indicating the inclusion status of the $n$-th prime, such as $1101$ for $\{2,3,7\}$ where 0 indicates the omission of 5. Specifying an infinite holdout set requires an infinite string of bits (since both the holdout and divisor set are infinite), and there is an uncountably infinite number of infinite binary strings. This, of course means that there are a great multitude of uncomputable infinite holdout problems since the set of computable functions is countable. However, there exist potentially meaningful computable subclasses, such as using remainders of primes in a particular modulus as a criterion for their inclusion in the holdout set. Section 4 provides an example of such a problem. To encode one of these problems as a unique integer, let $N$ be the modulus, $c$ be the coefficient from $cn +1$ and let $B$ be a binary string of length $N$. Each bit in $B$ will be set to 1 if primes congruent to that remainder mod $N$ are in the holdout set. We can then convert $B$ into its integer equivalent, and store the problem as $2^{B} \times 3^N \times 5^c$. There are probably more compact ways of doing this, but this is sufficient to produce unique integers. 

\subsection{The Finite Holdout Conjecture and its Implications}

Divisor and holdout problems are relatively simple to implement naively in code, and the subsequent optimization of this code can be a fun exercise. One possible experiment that the readers may choose to conduct is to take the first few primes, say $\{ 2, 3 ... 17\}$ pick one to be $c$ and try every combination of the others as a holdout set. The result is a little peculiar - none of them seem to diverge. Considering that this range contains numerous examples of divergent divisor problems this is actually shocking. We can even construct divergent divisor problems by hand. Consider $C_{3,\{7,11\}}$, for each prime in the divisor set, we may then try to solve the congruence $3x + 1\equiv x \pmod{p}$. the solution is $\lfloor p/2 \rfloor$ for any prime larger than 3. In the case of this example, the answers are 3 and 5. We may then apply the Chinese remainder theorem to find the remainder modulo 77 that satisfies the smaller congruences, this gives 38. If we start iterating $C_{3,\{7,11\}}$ from 38, it must diverge, since the remainder mod 7 and 11 will persistently be 3 and 5. This apparent difference in the divergence behaviors of divisor and holdout problems leads us to pose the finite holdout conjecture (FHC): that no finite holdout problem can have a divergent trajectory\\*

If the finite holdout conjecture is true, it would have a number of implications. Recall from the introduction that disproving the Collatz conjecture with a loop is "easy": simply run a loop-finding algorithm, and it will find the loop in some finite amount of time. Divergence on the other hand requires a proof, and that always carries the risk that the statement is unprovable (even if the trajectory does indeed diverge). In that case if the problem is also loop-free, then it is actually unfalsifiable. But if we have FHC, then all holdout problems that are false are false because they contain cycles, which are finite, and if they exist can always be shown to be there in a finite number of steps by performing the loop. It could be that some holdout problem is true but unverifiable, but with FHC we know that all false holdout problems can be discovered.\\*

A notable property of holdout problems is that they are not only the "opposite" of divisor problems. They can also serve as divisor problem approximators. Consider the holdout problem $H_{3,\{3,5,7,11...101\}}$ which includes all the primes less than or equal to 101, excluding 2, in its holdout set. If we were to perform a single iteration of this holdout problem on any number less than 34, it will perfectly match the iteration of the Collatz function. This is because taking the $3n+1$ transition will produce a number less than 101, and for such a small number, preserving all the primes that are not 2 is equivalent to preserving just the holdout set. Problems like this one, with a large enough holdout set can follow the Collatz function for arbitrarily large ranges of integers.\\*

By a similar line of reasoning, we can show that if we have FHC, and take the union of the sets of prime divisors of the (odd) iterates of a divergent trajectory of the Collatz function, the resulting set must be infinite. If this were not the case, then we could use that set as a holdout set for a finite holdout problem. That holdout problem would behave identically to the Collatz function for that trajectory since every prime factor that appears will be a member of the holdout set. That trajectory, however, is divergent, which FHC disallows.\\*

Investigating the Collatz conjecture as the limit of a series of holdout problems, as we have seen above, may well prove to be an important lens from which to view the problem. Determining the truth of the finite holdout conjecture is of critical importance, not just to understand the Collatz conjecture, but also to understand holdout problems simply for their own sake. But perhaps the even greater task lies in addressing what seems to be the lack of understanding about how adding 1 affects the prime factorization of a number. This is one of the key differentiators between divisor problems Minsky machines, and if it could be shown that adding 1 had a random enough effect, the possibility of divisor problems performing computation could be ruled out.

\bibliographystyle{IEEEbib}
\bibliography{refs}
\end{document}